

\baselineskip=14pt
\parskip=10pt
\def\halmos{\hbox{\vrule height0.15cm width0.01cm\vbox{\hrule height
  0.01cm width0.2cm \vskip0.15cm \hrule height 0.01cm width0.2cm}\vrule
  height0.15cm width 0.01cm}}
\def\halmos{\hbox{\vrule height 1ex width 1ex}}
\font\eightrm=cmr8 

\magnification=\magstephalf

\def\M{{\cal M}}

\def\1{{\overline{1}}}
\def\2{{\overline{2}}}
\parindent=0pt
\overfullrule=0in

\def\frac#1#2{{#1 \over #2}}
\centerline
{\bf A Simple Re-Derivation of Onsager's Solution of the 2D Ising Model using Experimental Mathematics }

\bigskip
\centerline
{\it Manuel KAUERS and Doron ZEILBERGER}

{\bf Abstract}: In this {\it case study}, we illustrate the great potential  of experimental mathematics
and symbolic computation, by rederiving, {\bf ab initio}, Onsager's celebrated solution of the two-dimensional
Ising model in zero magnetic field. Onsager's derivation is extremely complicated and ad hoc, as are all the
subsequent proofs. Unlike Onsager's, our derivation is not rigorous, yet it 
is {\it absolutely certain} (even if Onsager did not do it before),
and should have been acceptable to physicists who do not share mathematicians' fanatical (and often misplaced) insistence on rigor.

{\bf Two Warm-Up Exercises}

{\bf Definition 1}: For an  $n_1 \times n_2$ matrix, $M=(m_{i,j})$, and any positive real numbers $x$ and $y$:
$$
weight(M)\,(x,y) \, := \, x^{\frac{1}{2} \, \left ( \sum_{i,j} m_{i,j} \, m_{i+1,j} \, + \,  m_{i,j} \, m_{i,j+1} \right )} \, 
\cdot \,  y^{\sum_{i,j} m_{i,j}} \quad .
$$

(We make the convention that if $i$ is $n_1$, then $i+1=1$, and  if $j=n_2$ then $j+1=1$.)

{\bf Definition 2}:
Let $\M(n_1,n_2)$ be the set of $n_1 \times n_2$ matrices whose entries are either $1$ or $-1$ (of course, there are
$2^{n_1\,n_2}$ such matrices). The Laurent polynomial $P_{n_1,n_2}(x,y)$ is defined as follows.
$$
P_{n_1,n_2}(x,y) \, := \, \sum_{M \in \M(n_1,n_2)} \, weight(M)\,(x,y) \quad .
$$

{\bf Definition 3}: For $x,y$ positive real numbers:
$$
f(x,y) \, := \, \lim_{n \rightarrow \infty } \frac{ \log \, P_{n,n}(x,y)}{n^2} \quad .
$$

{\bf Exercise 1}: Find an explicit, {\bf closed-form}, expression for $f(x,y)$.

{\bf Definition 1a}: For an  $n_1 \times n_2 \times n_3$ three-dimensional array, $M=(m_{i,j,k})$, and a positive real number $x$,
$$
weight'(M)\,(x) \, := \, x^{\frac{1}{2} \, \left ( \, \sum_{i,j,k} m_{i,j,k} \, m_{i+1,j,k} \, + \,  m_{i,j,k} \, m_{i,j+1,k} \, + \, m_{i,j,k} \, m_{i,j,k+1}
\, \right) } \quad .
$$

{\bf Definition 2a}:
Let $\M(n_1,n_2,n_3)$ be the set of $n_1 \times n_2 \times n_3$ three-dimensional arrays, whose entries are either $1$ or $-1$ (of course, there are
$2^{n_1\,n_2\,n_3}$ such arrays), define the Laurent polynomial in $x$, by
$$
Q_{n_1,n_2,n_3}(x) \, := \, \sum_{M \in \M(n_1,n_2,n_3)} \, weight'(M)\,(x) \quad .
$$

{\bf Definition 3a}: For $x$, a positive real number,
$$
g(x) \, := \, \lim_{n \rightarrow \infty} \frac{ \log \, Q_{n,n,n}(x)}{n^3} \quad .
$$

{\bf Exercise 1a}: Find an explicit, {\bf closed-form}, expression for $g(x)$.

We hope, dear readers, that you will spend {\it some} time trying to solve these two exercises, but please do not
spend too much time! These have been open for almost eighty years, and in spite of many attempts by the best minds
in mathematical physics, are still wide open.

Exercise 1 is called ``solving the two-dimensional Ising model with magnetic field", while Exercise 1a is called
``solving the three-dimensional Ising model in zero magnetic field". Let us quote Ken Wilson, who got the
Physics Nobel prize in 1982 for seminal ({\it non-rigorous}!) work on these two `exercises' (without actually solving them!).

\qquad {\it ``When I entered graduate school I had carried out the instructions given to me by my father 
{\eightrm [notable chemist E. Bright Wilson, who co-authored, with Linus Pauling, the classic Introduction to Quantum Mechanics]}
and had knocked on both Murray Gell-Mann's and Feynman's doors and asked them what they were currently doing. Murray wrote
down the partition function for the three-dimensional Ising model and said that it would be nice if I could solve it. Feynman's answer was
`nothing' . ''}  \quad [Quoted in Julia Yeomans' wonderful book [Y], p. 35 .]

{\bf Onsager's Solution}

In 1944, Lars Onsager famously derived, and {\it rigorously} proved, the special case of Exercise 1, when $y=1$.

{\bf Onsager's Explicit Formula For the Zero-Field 2D Ising Model}: Let
$$
G(z) \, := \, -\frac{1}{4} \, \sum_{r=1}^{\infty} \, {{2r} \choose {r}}^2 \, \frac{z^{2r}}{r} \quad,
$$
then
$$
f(x,1) \, = \, \ln (x+x^{-1}) \, + \, G \left (\, \frac{x-x^{-1}}{(x +x^{-1})^2} \right ) \quad .
$$

Onsager's proof [O], and all the subsequent proofs, are very complicated. We will soon show how this formula could have
been {\it naturally} derived, way back in 1941, if they had the software and hardware that we have today (and even, probably, thirty years ago).

Unlike Onsager's derivation, that is fully rigorous, ours is not. 
So from a strictly (currently mainstream) mathematical viewpoint, it would have been considered `only' a conjecture,
were it done before Onsager's rigorous derivation.
But this conjecture would have been {\bf so} plausible that it would
have been whole-heartedly accepted by the theoretical physics community.

{\bf What is an ``Explicit" Answer}

From now on we will write $f(x)$ instead of $f(x,1)$, and $P_{n_1,n_2}(x)$ instead of $P_{n_1,n_2}(x,1)$.

In some sense Onsger's solution is disappointing and not really ``explicit", since it involves an infinite series,
that entails taking a limit. The {\it definition} of the function $f(x)$ also involves taking a limit
(namely of $\frac{\log(P_{n,n}(x))}{n^2}$ as $n \rightarrow \infty$). Why is the former limit better than the latter?

Indeed, the notion of ``explicit", or ``closed form" is vague and cultural. In ancient Greece a geometrical construction was
acceptable  only if it used ruler and compass. In algebra, for a long time, a solution was acceptable
only if it could be expressed in terms of the four elementary operations and root extractions. 
In enumerative cominatorics, a solution was (and sometimes still is) considered closed form only if it
is a product and/or quotient of factorials, and there are many other examples.

In a famous {\it position paper} [W], Herb Wilf tackled this problem in combinatorics. He was inspired to write
it when he was asked to referee a paper containing a ``formula" for a certain quantity. It turned out that
computing the quantity via the formula took much longer than using the definition. Inspired by the---at the time---new paradigm of ``computational complexity",
he suggested that an ``answer" is an efficient algorithm to compute the
quantity in question.

How would we compute $f(x)$, using the definition for a specific, `numeric', $x$?
We can, in principle, compute the sequence of Laurent polynomials $P_{n,n}(x)$ directly, for, say, $n \leq 30$,
and get the finite sequence of numbers $\{ \log P_{n,n}(x)/n^2 \}_{n=1}^{30}$, and see whether they get closer-and-closer, and estimate
the limit. Alas, computing $P_{n,n}(x)$ by brute force involves adding up $2^{n^2}$ terms, each of which take
$O(n^2)$ operations to compute. This is hopeless! Also, to be fully rigorous, one has to be able to find {\it a priori}
bounds for the error, and for each $\epsilon$ find (rigorously) an $n_{\epsilon}$ such that 
$|f(x)- \log(\, P_{n,n}(x) \, )/n^2|<\epsilon$ for $n \geq n_{\epsilon}$. This is truly hopeless.

On the other hand, using elementary calculus, Onsager's solution enables us to compute $f(x)$, very fast, to any 
desired accuracy.

More importantly, physicist do not really care about the explicit form of $f(x)$ (or more generally,
the still wide open  $f(x,y)$, and $g(x)$), they want to know the {\it exact} location of the {\bf singularities},
(critical points) that describes at what value of $x$ (and hence at what temperature) a  {\it phase transition} occurs, e.g. at
what temperature water boils or freezes. Even more importantly, they care about the
{\bf nature} of the singularities, in other words, {\it how} water boils 
rather than at what temperature (that depends, e.g. on pressure).
From Onsager's solution, one can easily find, using Calculus~I,
the location, and nature, of the singularity of $G(z)$, and hence of $f(x)$.
It is impossible to extract this information directly from the definition.

This motivation may be interesting, but it is irrelevant to us. All we want is to answer exercise 1 in the special case $y=1$,
with as little effort as possible, and making full use of the computer.
We only require elementary calculus and very elementary matrix algebra.
We don't even use eigenvalues!

{\bf Recommended Reading}

Even though it is {\it irrelevant} to our story,
for those readers who {\bf do} wish to know the context and background, we strongly recommend Barry Cipra's [C]
very lucid and very engaging introduction to the Ising model.  We also recommend the excellent books [T] and [Y].

{\bf Symbol-Crunching}

Of course, it would be nice to find an {\it expression} for $f(x)$ in terms of the {\bf symbol} $x$.
Computing $P_{n,n}(x)$ for any specific $n$ is a finite (albeit {\bf huge}) computation, involving summing
$2^{n^2}$ monomials, so we can't go very far. But, let's assume that we live in an ideal world,
or that quantum computing became a reality, then computing $P_{n_1,n_2}(x)$, and in particular, $P_{n,n}(x)$,
being finite, is always possible. The first, very natural, step, already proposed in 1941, that was motivated by the combinatorial approach 
(see later, and [T], Ch.6, Eq. 1.9, where we replace $x^2$ by $x$) is to write
$$
P_{n_1,n_2}(x)=  \frac{(x+2+x^{-1})^{n_1\,n_2}}{2^{n_1\,n_2}} Z_{n_1,n_2}(w) \quad, \quad{\rm where} \quad w=\frac{x-1}{x+1} \quad .
$$
It follows from a simple combinatorial argument that $Z_{n_1,n_2}(w)$ is
a {\it polynomial} in $w$, of degree $n_1\,n_2$.

Taking logarithms, and dividing by $n_1\,n_2$, we get
$$
\frac{\log P_{n_1,n_2}(x)}{n_1\,n_2} = -\log 2 + \log(x^{-1}+2+x) +\frac{\log Z_{n_1,n_2}(w)}{n_1\,n_2} \quad .
$$
Using the fact (do it!) that $x^{-1}+2+x=\frac{4}{1-w^2}$ we get that
$$
f(x)= \log 2 \, - \, \log(1-w^2) + \lim_{n \rightarrow \infty}\frac{\log Z_{n,n}(w)}{n^2} \quad{\rm where}\quad w=\frac{x-1}{x+1}.
$$
So from now, all we need is to find
$$
F(w):= \lim_{n \rightarrow \infty} \frac{\log Z_{n,n}(w)}{n^2} \quad .
$$

Now, it turns out (and it follows from elementary considerations) that the sequence  $\frac{\log Z_{n,n}(w)}{n^2}$ converges
in the sense of `formal power series'. More precisely, for any positive integer, $r$, the coefficient of $w^r$ in
$F(w)$ (our object of desire) coincides with that of  $\frac{\log Z_{n,n}(w)}{n^2}$ as soon as $n>r$. So a natural
{\it experimental mathematics} approach would be to try and find as many Taylor coefficients of $F(w)$ as our computer
would allow and look for a {\it pattern} that would enable us to conjecture a closed-form expression for
the Taylor coefficients of $F(w)$, thereby determining $F(w)$ and hence $f(x)$.

In an ideal world, with an indefinitely large computer, this very naive approach would have succeeded. Alas, as it turned out, we would
have needed to compute $P_{n,n}(x)$  for $n=96$, and since $2^{96^2}$ is such a big number, this
very naive {\bf brute force} approach is doomed to failure in our tiny universe.

{\bf Using Transfer Matrices}

A much more efficient approach to computing the Laurent polynomials $P_{n_1,n_2}(x)$ (and hence the polynomials $Z_{n_1,n_2}(w)$), 
was suggested in the seminal paper of Kramers and Wannier [KM]. That was also Onsager's starting point. It is easy to see
(see [T], p. 118) that for each $n_1$, there are easily computed $2^{n_1}$ by $2^{n_1}$ matrices, let's call them $A_{n_1}(x)$ such that
\def\Trace{{\rm Trace}}%
$$
P_{n_1,n_2}(x)= \Trace\, A_{n_1}(x)^{n_2} \quad .
$$
With today's computers, it is possible to compute these for $n_1 \leq 12$ and as large as $n_2$ as desired.

But once again, one can (still) not go very far.

In 1941, B.L. van der Waerden suggested an ingenious (very elementary!) combinatorial approach, described beautifully in Barry Cipra's article [C]
(see also Chapters 6 of [T] and [Y] for nice accounts).
He observed that the coefficients of $w$ in the polynomial $Z_{n_1,n_2}(w)$ have a nice combinatorial interpretation.
Putting $N=n_1\,n_2$, it turned out (and is very easy to see, see [T]) that for any positive integer $r$,
the coefficient of $w^r$ in $Z_{n_1,n_2}(w)$, let's call it $p_r$, is the number of `lattice polygons' with 
$r$ edges that can lie in an $n_1$ by $n_2$ `torodial rectangle', i.e. the set $\{0, \dots, n_1\} \times \{0, \dots, n_2\}$ with $0$ identified
with $n_1$ and $n_2$ respectively. A lattice polygon is a collection of edges such that every participating vertex has an {\it even} number
($0$, $2$, or $4$) of neighbors. It follows in particular that $p_r$ is zero if $r$ is odd.

It also follows from elementary combinatorial considerations that for $n_1>r,n_2>r$,
the coefficient $p_{r}$ is a certain {\it polynomial} in $N$ ([T], p. 150, Eq. (1.17)), and hence may be written $p_r(N)$, and we can write:
$$
p_{r}(N)= N a_r^{(1)} + N^2 a_r^{(2)} + \dots +  N^m a_r^{(m)} \quad .
$$
Now it also follows from elementary considerations, already known in 1941, that once you take the log, divide by $N=n_1\,n_2$ and
take the limit, {\it only} the coefficients of $N$ in these `Ising polynomials' survive, and that
$$
F(w)=\lim_{n \rightarrow \infty} \frac{\log ( \, Z_{n,n}(w) \,)}{n^2} \,= \,  \sum_{r=0}^{\infty} a_r^{(1)} w^r \quad .
$$

It remains to compute as many Ising polynomials, $p_r(N)$, as our computers will allow us, extract the coefficients $a_r^{(1)}$ of $N$,
and hope to detect a {\it pattern}, to enable us to conjecture the general coefficient of $F(w)$,
and hence know $f(x)$.

{\bf How to compute the Combinatorial Ising Polynomials?}

The first thing that comes to mind, and works well for small $r$ is to actually look for  the kind of lattice
polygons that can show up, but as $r$ gets larger, this gets out of hand. 
Rather than do the intricate combinatorics, we
use the fact that  $P_{n_1,n_2}(x)=\Trace \, A_{n_1}(x)^{n_2}$, from which we can compute
$Z_{n_1,n_2}(w)$ for $n_1\leq 12$ (say) and $n_2$ as large as desired.
For each individual coefficient of $w^r$ ($r$ even),
we output it for sufficiently many specific $n_1$ and $n_2$, 
and then using {\it undetermined coefficients} or interpolation
we fit them into a polynomial (whose degree we know beforehand).
In fact, it is possible to get $p_{2r}(N)$ by looking at $n_1=r-2, n_2>r$, by excluding obvious polygons that belong to the
$(r-2) \times n_2$ torodial rectangle but are impossible for a larger rectangle.

{\bf The Ising Polynomials}

By using this  very naive approach (only using matrix multiplication and then taking the trace) our beloved computers came up with the
following first 10 Ising polynomials (we were able to find quite a few more, but as we will soon see, the
first ten polynomials suffice).
$$
p_2(N)= 0, \quad
p_4(N)= N, \quad
p_6(N)= 2\,N, \quad
p_8(N)=\frac{1}{2} \,N \left( 9+N \right), \quad
p_{10}(N)= \,N \left( 6+N \right),
$$
$$
p_{12}(N)= \frac{1}{6}\,N \left( 7+N \right)  \left( 32+N \right), \quad
p_{14}(N)=N \left( 130+21\,N+{N}^{2} \right),
$$
$$
p_{16}(N)= \frac{1}{24}\,N \left( 11766+1715\,N+102\,{N}^{2}+{N}^{3} \right), \quad
p_{18}(N)= \frac{1}{3} \,N \left( 5876+776\,N+49\,{N}^{2}+{N}^{3} \right),
$$
$$
p_{20}(N)= {\frac {1}{120}}\,N \left( 980904+118830\,N+7415\,{N}^{2}+210\,{N}^{3}+{N}^{4} \right) \quad .
$$
Extracting the coefficients of $N$, we get
$$
0, 1, 2, \frac{9}{2}, 12, \frac{112}{3} , 130, \frac{1961}{4}, \frac{5876}{3}, \frac{40871}{5} \quad .
$$
Hence $F(w)$ starts with
$$
F(w)=
{w}^{4}+2\,{w}^{6}+\frac{9}{2}\,{w}^{8}+12\,{w}^{10}+{\frac {112}{3}}\,{w}^{12}+130\,{w}^{14}+{\frac {1961}{4}}\,{w}^{16}+{\frac {5876}{3}}\,{w}^{18}+{\frac {40871}
{5}}\,{w}^{20} + \cdots\quad .
$$
However, these ten terms (and even forty of them) do not suffice to guess a pattern.

{\bf Duality Saves the Day}

Way back in 1941, in the seminal paper of Kramers and Wannier, that we have already mentioned, they discovered the
{\it duality relation} (see [C] for a lucid explanation)
$$
f \left( \, \frac{x+1}{x-1} \, \right) = f(x) - \log \left( \, \frac{x- x^{-1}}{2} \, \right ) \quad .
$$
Letting
$$
x^{*}=\frac{x+1}{x-1}  \quad ,
$$
the duality relation can be written as
$$
f(x^*) = f(x) - \log \left ( \, \frac{x- x^{-1}}{2} \, \right ) \quad ,
$$
or in a more symmetric form
$$
 f(x) - \log(x+x^{-1}) \, = \,  f(x^*) - \log(x^*+(x^*)^{-1}) \quad .
$$
It follows that a more natural, and hopefully user-friendly, function to consider is
$$
\bar{f}(x) :=  f(x) - \log(x+x^{-1}) \quad,
$$
and we have that  $\bar{f}(x)$ is unchanged under the involution $x \leftrightarrow x^*$,
$$
\bar{f}(x^*)=\bar{f}(x) \quad .
$$
It is natural to change from the variable $w$ to one that is invariant under the change
$x\leftrightarrow x^*$. There are many possibilities. Obviously, in order to ensure the
invariance, we can set $z=R(x,x^*)$ for any {\it symmetric} rational function~$R$.
We only need to ensure that when $w$ is expressed as a series in~$z$, this series has
positive order, so that we are allowed to substitute it into~$F(w)$.
Since $F(w)$ has only even exponents, we may also prefer that the series $w=w(z)$ has
only odd exponents in~$z$, so that the substitution does not introduce odd exponents
into~$F(w)$.

If we try a {\it template} (`ansatz')
$$
  z = \frac{a_{0,0} + a_{1,0} (x + x^*) + a_{0,1} x x^* + a_{2,0} (x+x^*)^2 + a_{1,1} (x+x^*)xx^* + a_{0,2} (xx^*)^2}{b_{0,0} + b_{1,0} (x + x^*) + b_{0,1} x x^* + b_{2,0} (x+x^*)^2 + b_{1,1} (x+x^*)xx^* + b_{0,2} (xx^*)^2}
$$
with undetermined coefficients $a_{i,j}$ and $b_{i,j}$, we get a system of polynomial equations that can be
easily solved using so-called Gr\"obner bases.
This gets translated  into an equation relating $z$ and $w$ by eliminating $x$,
using the fact that $x=\frac{1+w}{1-w}$. The (computer-generated) result is an equation of the form
$$
 (\dots) + (\dots)w + (\dots)w^2 + (\dots)w^3 + (\dots)w^4
+  (\dots)z + (\dots)wz + (\dots)w^2z + (\dots)w^3z + (\dots)w^4z = 0\quad ,
$$
where the dots stand for certain linear combinations of the undetermined coefficients
which we suppress here because of their size.
In order to ensure that the solution for $w$ of this equation is a series in $z$ with odd exponents only,
it suffices to force the coefficients of all terms $w^iz^j$ with $i+j$ even to zero.
This gives a linear system whose solution brings the equation down to
$$
  (w-1)w(w+1)(a_{0,0}+a_{0,1}+a_{0,2}) + (1+w^2)^2 z (b_{0,0}-b_{1,0}+b_{2,0})=0. 
$$
This suggests the choice
$$
z \, = \, \frac{cw(1-w^2)}{(1+w^2)^2},
\quad{\rm or}\quad
w \, = \frac{z}{c} + \frac{3z^3}{c^3} + \frac{22z^5}{c^5} + \frac{211z^7}{c^7} + \frac{2306z^9}{c^9} + \cdots,
$$
for some nonzero constant~$c$.
The value of $c$ is not important. We take $c=2$ in order to cancel the term $\log 2$ below. 

Let $\bar{f}(x)$, in terms of $w$ be written $\bar{F}(w)$, then (since $x+x^{-1}=\frac{2(1+w^2)}{(1-w^2)}$; note that $x=\frac{1+w}{1-w}$)
$$
\bar{F}(w) :=  f(x) - \log(x+x^{-1}) =-\log(1-w^2)+F(w) +\log 2 -\log \left (\,  \frac{2(1+w^2)}{1-w^2} \, \right )
$$
$$
= \, -\log(1+w^2) \, + \, \sum_{r=0}^{\infty} a_r^{(1)} w^r  \quad ,
$$
giving
$$
\bar{F}(w) \, = \,
{w}^{4}+2\,{w}^{6}+{\frac{9}{2}}\,{w}^{8}+12\,{w}^{10}+{\frac {112}{3}}\,{w}^{12}+130\,{w}^{14}+{\frac {1961}{4}}\,{w}^{16}+{\frac {5876}{3}}\,{w}^{18}+{\frac {40871}{5}}\,{w}^{20} + O(w^{22}) \quad .
$$
Changing the variable to $z$, and renaming $\bar{F}(w)$ to $G(z)$, we get
$$
G(z)=-\frac{1}{4}\,{z}^{2}-{\frac {9}{32}}\,{z}^{4}-{\frac {25}{48}}\,{z}^{6}-{\frac {1225}{1024}}\,{z}^{8}-
{\frac {3969}{1280}}\,{z}^{10}-{\frac {17787}{2048}}\,{z}^{12
}
$$
$$
-{\frac {184041}{7168}}\,{z}^{14}-{\frac {41409225}{524288}}\,{z}^{16}-{\frac {147744025}{589824}}\,{z}^{18}-{\frac {2133423721}{2621440}}\,{z}^{20} +O(z^{22}) 
\quad .
$$
The first ten terms of the sequence of coefficients, let's call them $\{b_{2r}\}_{r=1}^{10}$ 
$$
-\frac{1}{4}, -\frac{9}{32}, - \frac{25}{48}, - \frac{1225}{1024}, -\frac{3969}{1280}, -\frac{17787}{2048}, -\frac{184041}{7168}, -\frac{41409225}{524288}, 
-\frac{147744025}{589824}, -\frac{2133423721}{2621440}, \dots
$$
factorizes nicely, and there is an obvious pattern. By fitting the sequence of {\bf ratios} $\{b_{2r+2}/b_{2r}\}_{r=1}^{9}$
into a rational function, the computer guesses
$$
\frac{b_{2r+2}}{b_{2r}} \, = \, \frac{r(2r+1)^2}{(r+1)^3} \quad ,
$$
that implies the closed-form expression, for the coefficients
$$
b_{2r}= -\frac{ {{2r} \choose {r}}^2}{r 4^{r+1}} \quad .
$$
Since we can (nowadays!) easily extend the sequence $b_{2r}$ up to (at least) sixteen terms, and this
`guess' indeed continued to hold, this makes it virtually certain that the guess is correct.
Combining everything, we derived, {\it ab initio}, by {\it pure guessing} (and very elementary and natural reasoning),
Onsager's formidable formula. \halmos 

{\bf What's next?}

Now that we have rediscovered Onsager's explicit formula for $f(x)=f(x,1)$,
a natural next step towards the general case $f(x,y)$ is to determine an
explicit expression for $m(x)=\frac d{dy} f(x,y)|_{y=1}$, i.e., the next term
in the Taylor series expansion of $f(x,y)$ with respect to~$y$ at $y=1$. Physicists call
this the {\it ``spontaneous magnetization''}.

Using transfer matrices, as before, it is easy to compute the first few
terms of $m(x)$ as a series in~$x$ (or $w$, or~$z$),
and we don't even need a computer to guess an explicit expression for them:
they all are zero. But that's just a part of the story.

Onsager observed that $m(x)$ is only zero for $x<1+\sqrt2$,
while for $x\geq1+\sqrt2$, it is equal to
$$
  \left(\frac{(x^2+1)^2 (x^2-2 x-1) (x^2+2 x-1)}{(x-1)^4(x+1)^4}\right)^{1/8}\quad.
$$
According to Thompson ([T], p.~135), this expression {\it ``was first derived by Onsager in
the middle of the 1940s, but in true Onsager fashion he has not to this day published
his derivation''.}

We don't know how he found this expression, but here is one way one could search for it, using
experimental mathematics.
For specific numbers $x,y$, we can compute numerical approximations of $f(x,y)$ using
the original definition (Def.~3 above).
For example, taking $f(x,y)\approx \log P_{n,n}(x,y)/n^2$ with $n\approx 20$ gives
several correct digits at a reasonable computational cost.
From the numerical estimates of $f(x,y)$ for various points $x,y$, we can obtain
numerical estimates for $m(x)$ and $m'(x)$, for various points~$x$.

The idea is to fit a differential equation against this numeric data.
Suppose we suspect a differential equation of the form
$$
  (a_0 + a_1 x + \cdots + a_{10}x^{10})m(x)
 +(b_0 + b_1 x + \cdots + b_{10}x^{10})m'(x) = 0,
$$
with unknown integer coefficients $a_i,b_i$ to be determined.
So for a specific point~$x$, the task is to find a so-called {\it integer relation}
of the real numbers $m(x),\dots, x^{10}m(x),m'(x),\dots,x^{10}m'(x)$.
There are well-known algorithms for finding such relations [FB, LLL].

In order to recover the relation from the values at a single point~$x$, we would need to
compute these values to a rather high precision, which is not an easy thing to do.
We can get along with less precision by using several evaluation points and searching
for a simultaneous integer relation of the numbers
$m(x),\dots, x^{10}m(x),m'(x),\dots,x^{10}m'(x)$, for several~$x$.
It turns out that by using enough evaluation points, we just need
about 6 decimal digits of accuracy of $m(x)$ and~$m'(x)$ for each of these points,
in order to establish a convincing guess. 
Unfortunately, this is a still bit more than what we were able to obtain
by a direct computation via transfer matrices.

{\bf Supporting Software:}

For Maple and C programs, as well as output files, please visit the web-page

{\tt http://sites.math.rutgers.edu/\~{}zeilberg/mamarim/mamarimhtml/onsager.html} \quad .

{\bf References}

[C]  Barry A. Cipra, {\it An Introduction to the Ising Model}, The American Mathematical Monthly {\bf 94} (1987),  937-959. \hfill\break
Available from {\tt http://www.yaroslavvb.com/papers/cipra-introduction.pdf} \quad  (accessed May 9, 2018).

[FB] H.R.P. Ferguson and D.H. Bailey, {\it A Polynomial Time, Numerically Stable Integer Relation Algorithm.} RNR Techn. Rept. RNR-91-032, Jul. 14, 1992. 

[KW] H.A. Kramers and G.H. Wannier, {\it Statistics of the Two-Dimensional Ferromagnet. Part I},
Physical Review {\bf 60} (1941), 252-262. \hfill\break
Available from {\tt  http://sites.math.rutgers.edu/\~{}zeilberg/akherim/KW1941.pdf} (accessed May 9, 2018).

[LLL] A.K. Lenstra, H.W. Lenstra, and L. Lovasz, {\it Factoring Polynomials with Rational Coefficients.}
Math. Ann. 261, 515-534, 1982.

[O] Lars Onsager, {\it Crystal Statistics, I. A two-dimensional model with an order-disorder transition}, Physical Reviews {\bf 65}(1944), 117-149. \hfill\break
Available from {\tt http://www.phys.ens.fr/\~{}langlois/Onsager43.pdf} (accessed May 9, 2018).

[T] Colin J. Thompson, {\it ``Mathematical Statistical Physics''}, MacMillan, 1972.

[W] Herbert S. Wilf,{\it  What is an Answer?}, The American Mathematical Monthly {\bf 89}, 289-292.

[Y] J.M. Yeomans, {\it ``Statistical Mechanics and Phase Transition''}, Oxford University Press, 1992.

\bigskip
\hrule
\bigskip
Manuel Kauers,  Institute of Algebra, J. Kepler University Linz, Linz, Austria. \hfill\break
{\tt manuel dot kauers at jku dot at} \quad .

\bigskip
Doron Zeilberger, Department of Mathematics, Rutgers University (New Brunswick), Hill Center-Busch Campus, 110 Frelinghuysen
Rd., Piscataway, NJ 08854-8019, USA. \hfill\break
Email: {\tt DoronZeil at gmail dot com}   \quad .
\bigskip
\hrule
\bigskip
Version of May 14, 2018

\end